\documentclass{article}

\usepackage{pstricks}
\usepackage{amsmath, amsthm, amssymb}
\usepackage[all]{xy}
\usepackage{xspace}
\usepackage{algorithm}
\usepackage{algpseudocode}

\newtheorem{Theorem}{Theorem}[section]

\newtheorem{lemma}{Lemma}[section]
\newtheorem{remark}{Remark}[section]

\newcommand{\comment}[1]{}

\newcommand{\edge}[2]{\ensuremath{\left\{#1,#2\right\}}}

\newcommand{\TwoGreedy}{{\sc 2GREEDY}}

\usepackage{amsmath,amssymb,amsthm,textcomp}

\usepackage{enumerate}

\definecolor{brown}{cmyk}{0, 0.72, 1, 0.45}
\definecolor{grey}{gray}{0.5}

\def\z{\zeta}

\allowdisplaybreaks
\def\HAM{{\sc extend-rotate}}
\def\cG{{\cal G}}

\def\hA{\hat{A}}
\def\hB{\hat{B}}
\def\hC{\hat{C}}
\def\hD{\hat{D}}

\def\bu{{\bf u}}

\def\s{\sigma}

\def\ex{\mathbb{E\/}}

\def\la{\lambda}
\def\la{\lambda}
\def\a{\alpha}

\def\part{\partial}
\def\tag#1 {\eqno(#1)}

\newcommand{\brac}[1]{\left( #1 \right)}
\newcommand\bfrac[2]{\left(\frac{#1}{#2}\right)}

\def\a{\alpha} \def\b{\beta} \def\d{\delta} \def\D{\Delta}
\def\e{\epsilon} \def\f{\phi}   \def\g{\gamma}
\def\G{\Gamma}  \def\k{\kappa} 
\def\z{\zeta} \def\th{\theta}   \def\l{\lambda}
 \def\m{\mu} \def\n{\nu} 
\def\r{\rho}  \def\s{\sigma}

\def\bd{{\bf d}}
\def\bv{{\bf v}}

\def\M{M^*}

\def\cE{{\cal E}}

\def\e{\varepsilon}
\def\n{\nu}
\newcommand{\proofstart}{{\bf Proof\hspace{2em}}}
\newcommand{\proofend}{\hspace*{\fill}\mbox{$\Box$}}
\newcommand{\set}[1]{\left\{#1\right\}}
\newcommand{\ignore}[1]{}

\def\hy{\hat{y}}
\def\hz{\hat{z}}
\def\hu{\hat{\bu}}

\newcommand{\beq}[1]{\begin{equation}\label{#1}}
\def\eeq{\end{equation}}
\def\hm{\hat{\m}}
\def\hl{\hat{\l}}
\def\hT{\hat{T}}

\def\whp{\text{w.h.p.}}
\def\2G{{\sc 2greedy}}

\def\Pr{\mathbb{P}}
\def\E{\mathbb{E}}

\def\cP{{\cal P}}
\def\qs{\text{q.s.}}

\def\wL{\Lambda_0}
\def\gc{G_{n,m}^{\d\geq 3}}
\def\Gc{{\cal G}_{n,m}^{\d\geq 3}}
\title{An almost linear time algorithm for finding Hamilton cycles in sparse random graphs with minimum degree 
at least three.}

\author{Alan Frieze\thanks{Research supported in part by NSF Grant CCF2013110}\ \ and 
Simi Haber\\Department of Mathematical Sciences,\\
Carnegie Mellon University,\\Pittsburgh PA15217.}

\begin{document}
\maketitle

\begin{abstract}
We describe an algorithm for finding Hamilton cycles in random graphs. Our model is the 
random graph $G=\gc$. In this model $G$ is drawn uniformly from graphs with vertex set $[n]$, $m$ edges and
minimum degree at least three. We focus on the case where $m=cn$ for constant $c$. 
If $c$ is sufficiently large then our algorithm runs in $O(n^{1+o(1)})$ time and succeeds \whp
\end{abstract}
\section{Introduction}
The threshold for the existence of Hamilton cycles in random graphs has been known very
precisely for some time, Koml\'os and Szemer\'edi \cite{KoSz}, Bollob\'as \cite{B1}, Ajtai, Koml\'os and 
Szemer\'edi \cite{AKS}.
Computationally, the Hamilton cycle problem is one of the original NP-complete problems described in the 
paper of Karp \cite{K1}.
On the other hand Angluin and Valiant \cite{AV} were the first to show that the Hamilton cycle problem 
could be solved efficiently
on random graphs. The algorithm in \cite{AV} is randomised and very fast, $O(n\log^2n)$
time, but requires $Kn\log n$ random edges for sufficiently large $K>0$. 
Bollob\'as, Fenner and Frieze \cite{BFF} gave a deterministic polynomial time algorithm that 
works \whp\ at the exact threshold for Hamiltonicity, it is shown to run in $O(n^{3+o(1)})$ time.

The challenge therefore is to find efficent algorithms for graphs with a linear number of edges. 
Here we have to make some extra assumptions
because a random graph with $cn$ edges is very unlikely to be Hamiltonian. It will have isolated vertices.
It is natural therefore to consider models of random graphs with a linear number of edges and minimum degree
$\d$ at least two. In fact minimum degree three is required to avoid the event of having three vertices of degree 
two having a 
common neighbor. For example, in the case of random $r$-regular graphs, $r=O(1)\geq 3$,
Robinoson and Wormald \cite{RW1}, \cite{RW2} settled the existence question and Frieze, Jerrum, Molloy, Robinson 
and Wormald \cite{FJMRW} gave a polynomial time algorithm for finding a Hamilton cycle.
The running time of this algorithm was not given explicitly, but it is certainly $\Omega(n^3)$.

We will work on a model where the assumption is that $\d\geq 3$ as opposed to all vertices having degree exactly
three. It is tempting to think that existence results for the regualr case $r=3$ will help. Unfortunately, this is
not true.
The random graph $\gc$ is uniformly sampled from the set $\Gc$ of graphs with vertex set $[n]$, $m$ edges and 
minimum degree at least three.

Frieze,  \cite{F1} gave an
$O(n^{3+o(1)})$ time algorithm for finding
large cycles in sparse random graphs and this can be adapted to find Hamilton
cycles in $G_{n,cn}^{\d\geq 3}$ in this time for sufficiently
large $c$. The paper \cite{F2} gives an algorithm that 
reduces this to $n^{1.5+o(1)}$ for sufficiently large $c$. The main aim of this paper
is to construct an almost linear time algorithm for this model.
\begin{Theorem}\label{th1}
If $c$ is sufficiently large then our algorithm finds a Hamilton cycle in $\gc$, $m=cn$,
and runs in $O(n^{1+o(1)})$ time and succeeds \whp
\end{Theorem}
\begin{remark}\label{rem-}
The $n^{o(1)}$ term here is $(\log n)^{O(\log\log n)}$ which is tantalisingly close to best possible(?) $\log^{O(1)}n$.
\end{remark}

\section{Outline of the paper}
The paper \cite{CFM} gave an efficient algorithm for finding the maximum matching in a sparse
random graph. Its approach was to (i) use the simple greedy algorithm of Karp and Sipser \cite{KS}
and then (ii) augment it to a maximum matching using alternating paths. In this paper we replace
the Karp-Sipser algorithm with the algorithm \2G\ that \whp\ finds a 2-matching in $G=\gc$ with $O(\log n)$ 
components and we replace alternating paths with extensions and rotations. (A 2-matching is a spanning
subgraph of maximum degree at most two).

In Section \ref{GS} we will describe our algorithm. We will describe it in two subsections.  
We will describe \2G\ for finding a good 2-matching $M$ in detail in Section \ref{2G}. In section \ref{erot} we will
describe an algorithm \HAM\ that uses extensions and rotations to convert $M$ into a Hamilton cycle.
In Section \ref{dist} we discuss some ``residual randomness'' left over by \2G.
In Section \ref{gc} we prove some structural properties of $\gc$. 
In Section \ref{BFS} we prove some properties relating the output of \2G\ to the execution of \HAM.
In Section \ref{ftp} we do a final calculation to finsih the proof. In Section \ref{diff} we point
to our difiiculties in proving $n\log^{O(1)}n$ and in Section \ref{FM} we make some final remarks.
\section{Algorithm}\label{GS}
As already stated, there are two phases to the algorithm. First we find a good 2-matching $M$ and then we
convert it to a Hamilton cycle. We look first at how we find $M$.
\subsection{Algorithm \2G}\label{2G}
We greedily and randomly choose edges to add to $M$. 
Edges of $M$ are deleted from the graph. We let $b(v)\in \set{0,1,2}$ denote the degree of $v$ in $M$. 
Once $b(v)=2$ its incident edges
are no longer considered for selection. The vertex itself is deleted from the graph. 
Thus the graph from which we select edges will shrink as the algorithm progresses.
We will use $\G$ to denote the current subgraph from which edges are to be selected.
When there are vertices $v$ of degree $2-b(v)$ (or less) in $\G$, we take care to choose an edge incident with such
a vertex. Our observation being that there is a maximum cardinality 2-matching of $\G$ that contains such an edge.

If every vertex $v$ of $\G$ had degree at least $3-b(v)$ then we choose an edge randomly from edges that are incident with 
vertices $v$ that have $b(v)=0$. In this way, we quickly arrive at a stage where every vertex of $\G$ 
has $b(v)=1$. At this point we use the algorithm of \cite{CFM} to find a (near) perfect matching $M^*$,
which we add to $M$ as our final 2-matching.

We describe \2G\ in enough detail 
to make some of its claimed properties meaningful.

We let
\begin{itemize}
\item $\m$ be the number of edges in $\G$,
\item $Y_k=\set{v\in [n]:d_\G(v)=k\text{ and }b(v)=0}$, $k=0,1,2$,
\item $Z_k=\set{v\in [n]:d_\G(v)=k\text{ and }b(v)=1}$, $k=0,1$,
\item $Y=\set{v\in [n]:d_\G(v)\geq 3\text{ and }b(v)=0}$,\qquad
\item $Z=\set{v\in [n]:d_\G(v)\geq2\text{ and }b(v)=1}$,\qquad
\item $M$ is the set of edges in the current 2-matching.
\end{itemize}
Note that $V(\G)=[n]\setminus (Y_0\cup Z_0)$ and that $b(v)\in\set{0,1}$ for $v\in V(\G)$.

We will assume that the input to our algorithm is an ordered sequence $\s_m=(e_1,e_2,\ldots,e_m)$ where $m=cn$. 
Here $E_m=\set{e_1,e_2,\ldots,e_m}$ are the edges of $\gc$ and $\s_m$ is a random ordering of $E_m$. 
Once these orderings
are given, the vertices and edges are processed in a deterministic fashion. Thus for example, if \2G\ requires a 
random edge with
some property, then it is required to take the first available edge in the given ordering.

We now give details of the steps of \vspace{.1in}

\noindent
{\bf Algorithm \2G:}\vspace{-.1in}
\begin{description}
\item[Step 1(a) $Y_1\neq\emptyset$]\ \\
Choose $v\in Y_1$. 
We choose $v$ by finding the first edge in the ordering $\s$ that contains a 
member of $Y_1$.
Suppose that its neighbour in $\G$ is
$w$. We delete the edge $(v,w)$
from $\G$ add $(v,w)$ to $M$ and move $v$ to $Z_0$.
\begin{enumerate}[(i)]
\item If $w$ is currently
in $Y$ then move it to $Z$.
If it is currently in $Y_1$ then move it to $Z_0$.
If it is currently in $Y_2$ then move it to $Z_1$. Call this re-assigning $w$.
\item If $b(w)=1$ then we move $w$ to $Z_0$ and
make the requisite changes
due to the loss of other edges incident with $w$. Call this {\em tidying up}.
\end{enumerate}

\item[Step 1(b): $Y_1=\emptyset$ and $Y_2\neq\emptyset$]\ \\
Choose $v\in Y_2$. We choose $v$ by finding the first edge in the ordering $\s$ that contains a 
member of $Y_2$. Suppose that its neighbours in $\G$ are
$w_1,w_2$.

We choose one of the neighbors at random, say $w_1$.
We move $v$ to $Z_{1}$.
We delete the edge $(v,w_1)$
from $\G$ and place it into $M$.
In addition,
\begin{enumerate}[(i)]
\item If $b(w_1)=0$ then put
$b(w_1)=1$ and add the edge $(v,w_1)$ to $M$. Re-assign $w_1$.
\item If $b(w_1)=1$ then we delete $w_1$ from $\G$.  Tidy up.
\end{enumerate}

\item[Step 1(c): $Y_2=\emptyset$ and $Z_1\neq\emptyset$]\ \\
Choose $v\in Z_1$. We choose $v$ by finding the first edge in the ordering $\s$ that contains a 
member of $Z_1$. Let $u$ be the other endpoint of the path $P$ of $M$
that contains $v$. Let $w$ be the unique neighbour of $v$ in $\G$. We delete $v$ from $\G$ and
add the edge $(v,w)$ to $M$. In addition there
are two cases.
\begin{enumerate}[(1)]
\item If $b(w)=0$ then we re-assign $w$.
\item If $b(w)=1$ then we delete vertex $w$ and tidy up. 
\end{enumerate}

\item[Step 2: $Y_1=Y_2=Z_1=\emptyset$ and $Y\neq\emptyset$]\ \\
Choose the first edge $(v,w)\in E(\G)$ in the order $\s_m$ incident with a vertex $v\in Y$.
We delete the edge $(v,w)$ from $\G$ and add it to $M$. We move $v$
from $Y$ to $Z$. There are two cases.
\begin{enumerate}[(i)]
\item If $b(w)=0$ then move $w$ from $Y$ to $Z$. 
\item If $b(w)=1$  then we delete vertex $w$ and tidy up. 
\end{enumerate}
\item[Step 3: $Y_1=Y_2=Z_1=Y=\emptyset$]\ \\
At this point $\G$ will be distributed as $G_{\n,\m}^{\d\geq 2}$
for some $\n,\m$ where $\m=O(\n)$.
As such, it contains a (near) perfect matching $\M$ \cite{FP} and it can
be found in $O(\n)$ expected time \cite{CFM}. This step comprises
\begin{enumerate}[Step 3a]
\item Apply the Karp-Sipser algorithm to $\G$. W.h.p. this results in the construction of a matching 
$\M_1$ that covers all but $\tilde{O}(\n^{1/5})$ vertices $U=\set{u_1,u_2,\ldots,u_\ell}$. 
\item Now find augmenting paths from $u_{2i-1}$ to $u_{2i}$ for $i\leq \ell/2$. This produces the matching $\M$.
\end{enumerate}

\end{description}
The output of \2G\ is set of edges $M\gets M\cup \M$.

\subsection{Extension-Rotation Algorithm}\label{erot}
We now describe an algorithm,
\HAM\ that \whp\ converts $M$ into a Hamilton cycle.
The main idea is that of a {\em rotation}. Given a path $P=(u_1,u_2,
\ldots,u_k)$ and an edge $e=(u_k,u_i)$ where $i\leq k-2$ we say that the path
$P'=(u_1,\ldots,u_i,u_k,u_{k-1},\ldots,u_{i+1})$ is obtained from $P$ by a rotation.
$u_1$ is the {\em fixed} endpoint of this rotation. We say that $e$ is the {\em inserted} edge.

Given a path $P$ with endpoints $a,b$ we define a {\em restricted rotation search}
$RRS(\n)$ as follows: Suppose that we have a path $P$ with endpoints $a,b$. We start
by doing a sequence of rotations with $a$ as the fixed endpoint. Furthermore
\begin{enumerate}[R1]
\item We do these rotations in ``breadth first manner'', described in detail in Section \ref{BFS}.
\item We stop this process when we have either (i) created $\n$ endpoints or (ii)
we have found a path $Q$ with an endpoint that has a neighbor $w$
outside of $Q$. The path $Q+w$ will be longer than $P$.
We say that we have found an
{\em extension}.
\end{enumerate}
Let $END(a)$ be the set of endpoints, other than $a$, produced by this procedure. 
Assuming that we did not find an extension and having constructed $END(a)$, we take each
$x\in END(a)$ in turn and starting with the path $P_x$
that we have found from $a$ to $x$, we carry out R1,R2 above with $x$ as the fixed
endpoint and either find an extension or create a set of $\n$ paths with
$x$ as one endpoint and the other endpoints comprising a set $END(x)$ of size $\n$.

Algorithm \HAM
\begin{enumerate}[Step ER1]
\item Let $K_1,K_2,\ldots,K_r$ be the components of $M$ where $|K_1|=\max\set{|K_j|:j\in [r]}$.
If $K_1$ is a path then we let $P_0=K_1$, otherwise we let $P_0=K_1\setminus\set{e}$ where $e$ is any edge of
$K_1$. 
\item Let $P$ be the component of the current 2-matching $M$ that contains $P_0$.
If $P$ is not a cycle, go directly to ER3.
If $P$ is a Hamilton cycle we are done. Otherwise there is a vertex $u\in P$ and
a vertex $v\notin P$ such that $f=(u,v)$ is an edge of $G$, assuming that $G$ is connected,
see Lemma \ref{conn1}. Let $Q$ be the component containing $v$. By deleting an edge of $P$ incident 
to $u$ and (possibly) and edge of $Q$ incident with $v$ and adding $f$
 we create a new path of length at least $|P|+1$ with vertex set equal to $V(P)\cup V(Q)$.
Rename this path $P$.
\item 
Carry out $RSS(\n)$ until either an extension is found or we have constructed
$\n$ endpoint sets.
\begin{description}
\item[Case a:]
We find an extension. Suppose that we construct a path $Q$ with endpoints $x,y$ such
that $y$ has a neighbour $z\notin Q$.
\begin{enumerate}[(i)]
\item If $z$ lies in a cycle $C$ then let $R$ be a path
obtained from $C$ by deleting one of the edges of $C$ incident with $z$. Let now
$P=x,Q,y,z,R$ and go to Step ER2.
\item If $z=u_j$ lies on a path $R=(u_1,u_2,\ldots,u_k)$ where the numbering is chosen so
    that $j\geq k/2$ then we let $P=x,Q,y,z,u_{j-1},\ldots,u_1$ and go to Step ER2.
\end{enumerate}
\item[Case b:] If there is no extension then we search for an edge $f=(p,q)$ such that
$p\in END(a)$ and $q\in END(p)$. If there is no such edge then the algorithm fails.
If there is such an edge we let $Q$ be the corresponding path from $p$ to $q$. We replace $P$ in our
2-matching by the cycle $Q+f$ and go to ER2.
\end{description}
\end{enumerate}

\subsection{Execution Time of the algorithm}
The expected running time of \2G\ is $O(n)$ and \whp\ it completes in $O(n)$ time with 
a 2-matching $M$ with at most $K_1\log n$ components for some constant $K_1>0$. This follows from the results of
\cite{CFM} and \cite{F2}.

To bound the execution time of \HAM\ we first observe that 
it follows from \cite{AV} that $RSS(\n)$ can be carried out in
$O(\n^2\log n)$ time. We will take
$$\n=n^{1/2+O(\e)}$$ 
where 
\beq{defeps}
\e=\frac{K(\log\log n)^2}{\log n}
\eeq
where $K$ is a sufficiently large positive constant and that $c$ is sufficiently large. 

We now bound the number of executions of $RSS(\n)$.
Each time we execute Step ER3, we either reduce the number of components by one or
we halve the size of one of the components not on the current path. So if the component
sizes of $M_0$ are $n_1,n_2,\ldots,n_\k$ then the number of executions of Step ER3 can be
bounded by
$$\k+\sum_{i=1}^\k \log_2n_i=O(\log^2n).$$
So the total execution time is \whp\ of order
$$n+(n^{1/2+O(\e)})^2\log^2n=O(n^{1+O(\e)}).$$
This clearly suffices for Theorem \ref{th1}.

We will now turn to discuss the probability that our algorithm succeeds after we have described \2G.
We remind the reader that the analysis assumes that $c$ is sufficently large.
\section{Residual Randomness}\label{dist}
Let $G$ be a graph with an ordering of its edges and consider a run
of \TwoGreedy on that graph. At every point of time each vertex is
in one of the sets $Y_0, Z_0$, $Y_1, Y_2, Z_1, Y$ and $Z$
as defined above. 

We let the set of vertices that were removed from the graph while in $Z$
be denoted by $R$. We call them ``regular vertices''. These vertices
are removed from $\Gamma$ in the execution of a {\bf{Step 1}} or {\bf{Step 2}} of
\2G\ and they are internal vertices of paths of $M$ at the start of Step 3. 

For a vertex $v$ let $t_v$ be the time at which \2G\ deletes $v$ from $\G$. 
Vertices $w$ that are not deleted before the start of Step 3 are given $t_w=\infty$. 
A vertex is {\em early} if $t_v\leq n^{1-\e}$ and {\em late} otherwise. 
An edge $e_i$ is {\em punctual} 
if $i\leq (1-\a)m$ and {\em tardy} otherwise, where $\a$ is a small positive constant.

When a vertex $v\in R$ gets matching degree two we take the incident non-matching
edge $e$ with the lowest index in $\s$ to be its $Z$-witness. The fact that $v\in Z$ just before this
happens implies that $e$ exists. We let
$W$ denote the set of $Z$-witnesses. We next defne two sets $R_0$ and $\wL$:

We let
$$R_0=\set{v\in R:\;v\text{ is early and the $Z$-witness of $v$ is punctual} }.$$
and 
$$\wL=\set{v:\text{$v$ has degree at least 
4 in $\G_{n^{1-\e}}$ }}.$$
We may now state and prove the main lemma of this section.

\begin{lemma}
In what follows $R_0,\wL$ are defined with respect to $G$ and an ordreing of its edges. 
Let $e = \edge{x}{y}$ be a tardy edge of $G$ where $x\in R_0$ and $y\in \wL$.  Let $G'$ be the graph
obtained from $G$ by deleting $e$.
Assume that running \2G\ on $G$ up until Step 3 gives a 2-matching $M$ and a witness set $W$ and running
\2G\ on $G'$ up until Step 3 gives $M'$ and $W'$. Then
$M = M'$ and $W = W'$.
\end{lemma}
\proofstart
We claim that up to time $t_x$, \2G\ will delete the same vertices and edges from $\G_t$ and $\G_t'$ and then 
delete $x$ from both. After this the two graphs will coincide and we are done. We do this by induction on $t$.
This is clearly true for $t=0$ and assume that $\G_t$ and $\G_t'$ differ only in $e$ and $t<t_x$.
Note that the induction hypothesis implies that the sets $Y_0,Y_1,\ldots,Z$ are the same in $\G_t,\G_t'$.
This is because deleting edge $e$ does not affect $x$'s status, because $e$ is after the $Z$-witness of $x$ 
in the order $\s$. It does not affect $y$'s status, because the degree of $y$ will be at least 3 after the deletion.

Because the sets $Y_0,Y_1,\ldots,Z$ are unchanged, the choice of step is the same in $\G_t$ and $\G_t'$.
The difference between the two graphs only affects the degrees of $x$ and $y$ and by construction, this is
never enough to change a choice of edge.
\proofend

\begin{remark}\label{rem1}
Suppose that $e_i=(v,w)$ and that (i) $v\in R_0$, (ii) $w\in\wL$ and 
(iii) $e_i$ is tardy. Then replacing $e_i$ by $(v',w')$ such that (i) $v'\in R_0$ and (ii) $w'\in \wL$ 
results in the the same output $M,W$.

The net effect of this is that if we condition on all edges except for the tardy 
edges between $R_0$ and $\wL$ then the unconditioned tardy $R_0:\wL$ edges are random. This is what
we mean by there being residual randomness.
\end{remark}
\section{Degree Sequence of $\gc$}\label{gc}
The degrees of the vertices in $G$
are distributed as truncated Poisson random variables $Po(\l;\geq 3)$, see for example
\cite{AFP}. More precisely we can generate the
degree sequence by taking random variables $Z_1,Z_2,\ldots,Z_n$ where
\beq{prob}
\Pr(Z_i=k)=\frac{\l^k}{k!f_3(\l)}\qquad for\ i=1,2,\ldots,n\ and\ k\geq 3,
\eeq
where $f_j(\l)=e^\l-\sum_{k=0}^{j-1}\frac{\l^k}{k!}$ for $j\geq 1$.

Then we condition on $Z_1+Z_2+\cdots Z_n=2m$. The resulting $Z_1,Z_2,\ldots,Z_n$ can be taken to 
have the same distribution as the degrees of $G$. This follows from Lemma 4 of \cite{AFP}.
If we choose $\l$ so that 
$$\E(Po(\l;\geq 3))=\frac{2m}{n}\text{ or }\frac{\l f_2(\l)}{f_3(\l)}=\frac{2m}{n}$$ 
then
the conditional probability, $\Pr(Z_1+Z_2+\cdots Z_\n=2m)=\Omega(1/\sqrt{n})$
and so we
will have to pay a factor of $O(\sqrt{n})$ for removing the
conditioning i.e. to use the simple inequality $\Pr(A\mid B)\leq \Pr(A)/\Pr(B)$.
(This factor $O(n^{1/2})$ can be removed but it will not be necessary to do this here).

The maximum degree $\D$ in $G$ is less than $\log n$ \qs\footnote{A sequence of events,
$\cE_n$ occurs {\em quite surely} (\qs) if $\Pr(\neg\cE_n)=o(n^{-C})$ for any $C>0$. } and equation (7) of \cite{AFP} 
enables us to claim that 
that if $\n_k,2\leq k\leq \log n$ is the number of vertices of
degree $k$ then \qs
\begin{equation}
\label{degconc}
\left|\n_k-\frac{n \l^ke^{-\l}}{k!f_3(\l)}\right|
\leq K_1\left(1+\sqrt{n \l^ke^{-\l}/(k!f_3(\l))}\right)\log n,\ 2\leq k\leq \log n.
\end{equation}
for some constant $K_1>0$.

In particular, this implies that if the degrees of the vertices in $G$ are $d_1,d_2,\ldots,d_n$ then \qs
\beq{s2}
\sum_{i=1}^n d_i(d_i-1)=O(n).
\eeq

Given the degree sequence we make our computations in
the configuration model, see Bollob\'as \cite{B2}.  
Let $\bd=(d_1,d_2,\ldots,d_n)$ be a sequence 
of non-negative integers with $2m=cn$.
Let $W=[2cn]$ be our set
of {\em points} and let $W_i=[d_1+\cdots+d_{i-1}+1,d_1+\cdots+d_{i}]$,
$i\in [n]$, partition $W$. The function $\f:W\to[n]$ is defined by
$w\in W_{\f(w)}$. Given a
pairing $F$ (i.e. a partition of $W$ into $m=cn$ pairs) we obtain a
(multi-)graph $G_F$ with vertex set $[n]$ and an edge $(\f(u),\f(v))$ for each
$\{u,v\}\in F$. Choosing a pairing $F$ uniformly at random from
among all possible pairings of the points of $W$ produces a random
(multi-)graph $G_F$.

This model is valuable because of the following easily proven fact:
Suppose $G\in \cG_{n,\bd}$, the set of (simple) graphs with vertex set $[n]$ and degree sequence \bd. Then
$$\Pr(G_F=G\mid G_F\mbox{ is simple})=\frac{1}{|\cG_{n,\bd}|}.$$
It follows that if $G$ is chosen randomly from $\cG_{n,\bd}$, then for any graph property $\cP$
\begin{equation}\label{simple}
\Pr(G\in \cP)\leq \frac{\Pr(G_F\in \cP)}{\Pr(G_F\mbox{ is simple})}.
\end{equation}
Furthermore, applying Lemmas 4.4 and 4.5 of McKay \cite{McKay} we see that if the degree sequence of $G$
satisfies \eqref{s2} then $\Pr(G_F\mbox{ is simple})=\Omega(1)$. In which case the configuration 
model can substitute for $\cG_{n,\bd}$ (and hence $\gc$) in dealing with events of probability $o(n^{-1/2})$. 

\begin{lemma}\label{lemgrow}
W.h.p. 
\begin{enumerate}[(a)]
\item $\gc$ contains no set $S\subseteq [n], 3\leq s=|S|\leq s_0=\frac{1}{5}\log_cn$ such that $S$ contains at least 
$s+1$ edges.
\item Let $W_1$ denote the set of vertices $v$ that are within 
distance $\ell_0=2\log\log n$ of a cycle of length at most $2\ell_0$ in $\gc$. Then w.h.p.
$|W_1|\leq n^{1/2}\log^{4\ell_0}n$.
\item W.h.p. there does not exist a connected subset of $K_c\log n
\leq s\leq n^{3/5}$ vertices that contain $s/10$ vertices of 
degree at most $30$. Here $K_c$ is some sufficiently large constant.
\end{enumerate}
\end{lemma}
\proofstart
(a) The expected number of sets $S$ containing $|S|+1$ edges can be bounded by 
\begin{align}
&O(n^{1/2})\sum_{s=3}^{s_0}\sum_{|S|=s}\sum_{D\geq 3s}\sum_{\substack{d_1+\cdots+d_s=D\\d_1,\ldots,d_s\geq 3}}
\prod_{i=1}^s\frac{\l^{d_i}}{f_3(\l)d_i!}\binom{D}{s+1}\bfrac{D}{cn-2s}^{s+1}\leq\label{line1}\\
&O(n^{1/2})\sum_{s=3}^{s_0}\sum_{|S|=s}\sum_{D\geq 3s}\bfrac{De}{s+1}^{s+1}\bfrac{D}{cn}^{s+1}
\frac{\l^Ds^D}{D!f_3(\l)^s}.\label{line2}
\end{align}
{\bf Explanation:} For \eqref{line1} we choose a set of size $s$ with vertices of degree $d_1,d_2,\ldots,d_s\geq 3$ and 
$d_1+\cdots+d_s=D$. The term $\prod_{i=1}^s\frac{\l^{d_i}}{f_3(\l)d_i!}$ (modulo $O(n^{1/2})$) accounts for the
probability of these degrees. We then choose $s+1$ configuration points and approximate the probablity that 
they are all paired with other points associated with $s$ by $\bfrac{D}{cn-2s}^{s+1}$.
We use $\sum_{d_1+\cdots +d_s=D}\prod_{i=1}^s\frac{1}{d_i!}=\frac{s^D}{D!}$ to get \eqref{line2}.

Continuing we observe that $(D/cs)^{2s+2}\leq \brac{1+\frac{3}{c}}^D$ for $D\geq 3s$. Thisis clearly true for
$D\leq cs$ and follows by induction on $D\geq cs$. 
Therefore,
$$\sum_{D\geq 3s}\frac{D^{2s+2}\l^Ds^D}{D!}\leq (cs)^{2s+2}e^{(\l+3) s}.$$
Plugging this into \eqref{line2} we get a bound of
\begin{align*}
&O(n^{1/2})\sum_{s=3}^{s_0}\sum_{|S|=s}\frac{ces^2}{n}\bfrac{ces^2e^{\l+3}}{n(s+1)f_3(\l)}^s\\
\leq &O\bfrac{cs^2}{n^{1/2}}\sum_{s=3}^{s_0}\bfrac{ne}{s}^s\bfrac{ces^2e^{\l+3}}{n(s+1)f_3(\l)}^s\\
\leq &O\bfrac{cs^2}{n^{1/2}}\sum_{s=3}^{s_0}\bfrac{ce^{\l+3}}{f_3(\l)}^s\\
\leq &O\bfrac{cs^2}{n^{1/2}}\sum_{s=3}^{s_0}(2ce^3)^s\\
=&o(1).
\end{align*}

(b) 
$$\E(|W_1|)\leq O(n^{1/2})\sum_{k\leq \ell_0,\ell\leq 2\ell_0}\binom{n}{k+\ell}
k\ell\bfrac{\D^2n}{2m-6\ell_0}^{k+\ell}\leq 2n^{1/2}\ell_0^2\log^{3\ell_0}n.$$
(We remind the reader that it is possible to remove the $O(n^{1/2})$ factor here. This would be worth doing
if we could reduce $\e$ to $O(\log\log n/\log n)$. This should become apparant in the proof of Lemma \ref{lem2},
equation \eqref{qazx}).

Part (b) follows from the Markov inequality. 

(c) For a fixed $s$, the probability such a set exists can be bounded by
$$O(n^{1/2})\sum_{|S|=s}\binom{s}{s/10}\sum_{D\geq 3s}\sum_{\substack{d_1+\cdots+d_s=D\\3\leq d_i,\,i\in [s]\\
d_i\leq 30,\,i\in[s/10]}}\prod_{i=1}^s\frac{\l^{d_i}}{d_i!f_3(\l)}\binom{D}{s-1}\bfrac{D}{cn}^{s-1}.$$
{\bf Explanation:} We choose a set $S$ and we let the degrees in $S$ be $d_1,d_2,\ldots,d_s$ where $D$
is the total degree. Since the induced subgraph is connected, it must contain a spanning tree. We weaken this to
it must contain $s-1$ edges. $\binom{D}{s-1}$ enumerates the lower numbered point of the edges and then $D^{s-1}$ 
enumerates the other possible endpoints and then $\bfrac{1}{cn-2s}^{s-1}=\frac{1+o(1)}{(cn)^{s-1}}$ 
bounds the probability the selected pairs exist.

We bound this by
\begin{align*}
&O(n^{1/2})\binom{n}{s}\binom{s}{s/10}\sum_{D\geq 3s}\binom{D}{s}\bfrac{D}{cn}^{s-1}f_3(\l)^{-s}
[x^D]\brac{\sum_{i=3}^{30}\frac{\l^ix^i}{i!}}^{s/10}f_3(\l x)^{9s/10}\\
&\leq O(n^{3/2})\bfrac{e}{s}^s(10e)^{s/10}\sum_{D\geq 3s}\bfrac{De}{s}^{s}\bfrac{D}{c}^{s-1}
\frac{1}{f_3(\l)^{s}(1+\xi)^{D}}
\brac{\sum_{i=3}^{30}\frac{\l^i(1+\xi)^i}{i!}}^{s/10}f_3(\l(1+\xi))^{9s/10}
\end{align*}
for any positive $\xi$.

Now if $\l$ is large then $f_3(\l)\geq e^{\l}/2$. Also, $f_3(\l(1+\xi))\leq e^{\l(1+\xi)}$. Furthermore,
$$\sum_{i=3}^{30}\frac{\l^i(1+\xi)^i}{i!}\leq 2\frac{\l^{30}e^{30\xi}}{30!}
\leq 2\bfrac{\l e^{1+\xi}}{30}^{30}.$$
We will take $\xi$ to be small but fixed. Then the bound becomes 

$$O(n^{3/2})(10e)^{s/10}\bfrac{e^2}{cs^2}^{s}\frac{2^{s/10}e^{9\l\xi s/10}}{e^{\l s/10}}
\brac{2\bfrac{\l e^{1+\xi}}{30}^{30}}^{s/10}
\sum_{D\geq 3s}\frac{D^{2s-1}}{(1+\xi)^D}.$$
We observe that if $u_D=\frac{D^{2s}}{(1+\xi)^D}$ then $u_{D+1}/u_D\leq \frac{1+3/D}{1+\xi}$ for $D\geq 9$.
So,
$$\sum_{D\geq 3s}\frac{D^{2s-1}}{(1+\xi)^D}\leq (3s)^{2s}\sum_{D\geq 3s}\frac{1}{(1+\xi)^D}\prod_{i=3s}^D\frac{i+3}{i}
\leq (3s)^{2s}\sum_{D\geq 3s}\frac{(D+3)^3}{(3s)^3(1+\xi)^D}\leq \frac{(3s)^{2s}}{\xi}.$$
For the last inequality we use the fact that $s$ is large and then $D^3\ll(1+\xi)^D$.

Continuing, we get a bound of 
\begin{align*}
&\leq O(n^{3/2})(10e)^{s/10}\bfrac{e^2}{cs^2}^{s}\frac{2^{s/10}e^{9\l\xi s/10}}{e^{\l s/10}}
\brac{2\bfrac{\l e^{1+\xi}}{30}^{30}}^{s/10}(3s)^{2s}\\
&=O(n^{3/2})\brac{\frac{9(40e)^{1/10}e^2}{ce^{\l(1-9\xi)/10}}\bfrac{\l e^{1+\xi}}{30}^{5}}^s\\
&=o(1)
\end{align*}
if we take $\xi=1/10$ and $c$ and hence $\l$ sufficiently large.

\proofend
\section{Finding a Hamilton cycle}\label{BFS} 
We assume that we have a path $P$ with endpoints $a,b$ and we do rotations with $a$ as the fixed endpoint
to try to find an extension. In the next section we show that if no extensions are found, then \whp\ 
we create sufficient endpoints other than $b$ on paths of length equal to $P$.
Throughout this description, we will assume that no extension is found i.e. all neighbors of endpoints
turn out to be vertices of $P$. We associate the search with something similar to an alternating tree of matching 
theory.

\subsection{Tree Growth}
In this section we describe our search for a longer path than $P$ 
using \HAM\ in terms of growing a tree structure, where each
vertex determines a new long path. We expose what happens \whp\ if  we fail to find an extension.
Let $A_0=\set{b}$ and let $B_0$ be the set of neighbors of $b$ on $P$, excluding $b$'s path neighbor.
We now define the sets $A_i,B_i,i=1,\ldots,$ and $C_i=\bigcup_{j\leq i}(A_j\cup B_j)$. Here every vertex $v$ in $A_i$
will be the endpoint of a path of the same length as $P$. It will be obtained from $P$ by exactly $i$ rotations
with $a$ as the fixed endpoint.
Fix $i\geq 0$ and let $A_i=\set{v_1,v_2,\ldots,v_k}$. We build $A_{i+1},B_{i+1}$ by examining 
$v_1,v_2,\ldots,v_k$ 
in this order.
Initially $A_{i+1}=B_{i+1}=\emptyset$ and we will add vertices as we process the vertices of $A_i$.
Fix $v=v_j$. We have a path $P_v$ with endpoints $a,v$. We consider two cases:

{\bf Case 1:} $|C_i|\leq L_0=\frac1{20}\log_cn$.\\
Let $N_v=\set{u_1,u_2,\ldots,u_d}$ be the neighbors of $v$, excluding 
its neighbor on $P_v$. We also exclude from $N_v$ those neighbors already in $B_{i+1}$ (as defined so far).
Let $w_l$ be the neighbor of $u_l$ on $P_v$ that lies between $u_l$ and $v$
for $l=1,2,\ldots,d$. Let $N'_v=\set{w_1,w_2,\ldots,w_d}$. We exclude from $N'_v$ those vertices already in $A_{i+1}$ 
(as defined so far).
We add $N_v$ to $B_{i+1}$ and $N'_v$ to $A_{i+1}$ and we add edges $(v,u_j)$ and $(u_j,w_j)$ to $T$.
The edge $(u_j,w_j)$ will be called a {\em lost} edge. 
Furthermore, we define $P_{w_j}=P_v+(v,u_j)-(u_j,w_j)$ and observe 
that $P_{w_j}$
has endpoints $b,w_j$.

{\bf Case 2:} $|C_i|>L_0$.\\
Now let $N_v=\set{u_1,u_2,\ldots,u_d}$ be its neighbors as above.
We now exclude from $N_v$ those neighbors already in $C_{i+1}$ (as defined so far) as well as those $u_j$ for which 
$w_j\in C_{i+1}$.
We define $N'_v$ and update $A_{i+1},B_{i+1},T$ with this restricted $N_v$. 

\ignore{We should be able to use the calculations from [7]. We can prove the existence/non-existence independent
of the ordering of the edges. We can do this until the tree is of moderate size. Then maybe we can show growth 
regardless of path? We can argue for the beginning growth, $O(\log_cn)$ because we know there are many endpoints, [11].}

We define the subgraph $T=T(P,b,k)$ as follows: It has vertex set $C_k$ plus the edges of the form $(v,u_j)$ and 
$(u_j,w_j)$ used above.
$T$ suggests a tree. It is usually a tree, but in rare cases it may be unicyclic. This follows from 
Lemma \ref{lemgrow}. When this happens, some $v\in A_i$ (Case 1) has a neighbor in $B_j,j\leq i$.

We see from this that w.h.p. $T$ defined prior to the lemma has at most 
one cycle. By construction, cycles of $T$ are contained 
in the first $i_0$ levels. If there are two cycles inside the first $i_0$ levels then there is a set
$S$ (consisting of the two cycles plus a path joining them) with at most $4i_0$ vertices and at least 
$|S|+1$ edges.

We argue next that w.h.p. $T$ can be assumed to grow to a certain size and we can control its rate of growth.
\begin{lemma}\label{lem1}
Let $\b$ be some small fixed positive constant.
If $c$ is sufficiently large, then for all paths $P$ and endpoints $b$ such that extension does not occur, w.h.p.
\begin{enumerate}[(a)]
\item There exists $k$ such that $|C_k|\geq L_0=\frac{1}{15}\log_cn$.
\item If $L_0\leq |C_k|\leq n^{.6}$ then $|A_{k+1}|\in [2(1-\b)c|C_k|,2(1+\b)c|C_k|]$,
{\em even if only punctual edges are used once $|A_k|$ reaches size at least $n^\e$}.
\item There exists $k_0=O(\log_cn)$ such that $|A_{k_0}|\in [(2c(1+\b))^{-1}n^{1/2+5\e},2c(1+\b)n^{1/2+5\e}]$.
\item Let $k_1=k_0-\ell_0$ where $\ell_0=2\log\log n$
and let $x\in A_{k_1}$. Let $S$ be the set of descendants of $x$ in $A_{k_0}$
and let $s=|S|$. Let $S_0=\set{y\in S:d(y)\geq30}$ and let $s_0=|S_0|$.
Then, where $W_1$ is as in Lemma \ref{lemgrow},
\begin{enumerate}[(i)]
\item $x\notin W_1$ and $s\geq (2c(1-\b))^{k_0-k_1}/4$ implies that $s_0\geq 99s/100$. 
\item $s\leq (2c(1+\b))^{k_0-k_1}\log n$.
\end{enumerate}
\end{enumerate}
\end{lemma}
\proofstart
(a) Lemma 2.1 of \cite{FP} proves the following: Suppose that $S$ is the set of endpoints that can be produced by 
considering all possible
sequences of rotations starting with some fixed path $P$ and keeping one endpoint fixed. Let $T$ be the set of 
external neighbors of $S$. Here $S\cap T=\emptyset$.
Then $|T|\leq 2|S|$ and $S\cup T$ conmtains strictly more than $|S\cup T|$ edges. The assumption here is that the 
graph involved has minimum degree at least
three. It follows from Lemma \ref{lemgrow} that $|S|\geq \frac{1}{15}\log_cn$. As a final check, if $|C_k|$ never 
reached $L_0$ then it would have explored
all possible sets of endpoints i.e. the breadth first search is no restriction.

(b) 
If the condition in (b) fails then the following structure appears: Let $\d=1$ if $T$ is not a tree and
0 otherwise. Let $EVEN(T)=\bigcup_{i=0}^kA_i$ and
$ODD(T)=\bigcup_{i=0}^kB_i$ where $k$ is the number of iterations involved in the construction of $T$.
Then with $|EVEN(T)|=l$ and $|N(EVEN(T))|=r$ we have
(i) $2(l-1)+\d$ edges of $T$ connecting \textit{EVEN(T)} to
\textit{Odd(T)} (ii) $r-l+1$ edges connecting \textit{EVEN(T)} to
$N(\textit{EVEN(T)}) \setminus \textit{Odd(T)}$ and (iii) none
of the $l(n-r-l)$ edges between
$\textit{EVEN(T)}$ and $V\setminus N(\textit{EVEN(T)})$
are present. 

Assume first that $T$ is actually a tree and that $l\leq n^\e$ so that the edges of $T$ need not be punctual.

Given the vertices of $T$ and $N(EVEN(T))$, the probability of the existence of a $T$
with $L_0\leq l\leq n^{.6}$ and $r\leq (1-\b)cl$ can be bounded by
\beq{expression}
O(\sqrt{n})\left(\frac{1}{2m - 2(l+r)}\right)^{l+r-1}
\;\sum_{\substack {d_i \geq 3,\,i\in[r+l-1] \\\sum_{i=1}^{l} d_i = r+l-1}}
\left(\prod_{i=1}^{l} \frac{\lambda^{d_i}d_i!}{d_i!f_3(\l)} \prod_{i=l+1}^{2l-1}
\frac{\lambda^{d_i}d_i(d_i-1)}{d_i!f_3(\l)}\prod_{i=2l}^{r+l-1}
\frac{\lambda^{d_i}d_i}{d_i!f_3(\l)}\right)
\eeq
{\bf Explanation:}
The probability that an edge exists between
vertices $u$ and $v$ of degrees $d_u$ and $d_v$, given the existence
of other edges in $T$, is at most $\frac{d_u' d_v'}{2m-2(l+r) +
  3}$ where $d_u'=d_u$ less the number of edges already assumed to be incident with $u$. 
Hence, given the degree sequence, the probability that $T$ exists is at most

\begin{align*}
&\left(\frac{1}{2m - 2(l+r) }\right)^{l+r-1} \prod_{i=1}^{l} d_i!
\prod_{i=l+1}^{2l-1} d_i(d_i-1) \prod_{i=2l}^{r+l-1}
d_i.
\end{align*}
(We dropped the +3 in $2m-2(l+r) +3$).

Here the first product corresponds to $EVEN(T)$, the second product corresponds to $ODD(T)$ and the 
final product corresponds to neighbours of $T$ (not in $T$).

We will implicitly use the fact that if $c$ is sufficiently large, then so is $\l$.

We now simplify the
expression \eqref{expression} obtained for the probability to

\begin{align}
&O(\sqrt{n})\left(\frac{1}{2m - 2(l+r) }\right)^{l+r-1}\times\nonumber\\
&\frac{\lambda^{2r+2l-3}}{f_3(\l)^{r+l-1}}
\sum_{\sum_{i=1}^{l} d_i = r+l-1}
\brac{\sum_{d_i \geq 3,\,i\in[r+l-1]}
\prod_{i=l+1}^{2l-1}
\frac{\lambda^{d_i-2}}{(d_i-2)!}
\prod_{i=2l}^{r+l-1}
\frac{\lambda^{d_i-1}}{(d_i-1)!}}   \nonumber\\
&\leq O(\sqrt{n})\brac{\frac{1}{2m - 2(l+r) }}^{l+r-1}\times\nonumber\\
&\frac{\lambda^{2r+2l-3}}{f_3(\lambda)^{r+l-1}}
\sum_{\substack{\sum_{i=1}^{l} d_i = r+l-1\\d_i\geq 3,i\in[l]}}\brac{
\brac{\prod_{i=l+1}^{2l-1}\sum_{d_i\geq 3}\frac{\lambda^{d_i-2}}{(d_i-2)!}}
\brac{\prod_{i=2l}^{r+l-1}\sum_{d_i\geq 3}\frac{\lambda^{d_i-1}}{(d_i-1)!}}} \nonumber\\
&\leq O(\sqrt{n})\left(\frac{1}{2m - 2(l+r) }\right)^{l+r-1}
\frac{\lambda^{2r+2l-3}}{f_3(\l)^{r+l-1}}\;\;\binom{r}{l}
f_1(\l)^{l-1}\;\;f_2(\l)^{r-l}\label{zx7}\\
&\leq O(\sqrt{n})\left(\frac{1}{2m - 2(l+r) }\right)^{l+r-1}
\frac{\lambda^{2r+2l-3}}{f_3(\l)^{r+l-1}}\left(\frac{er}{l}\right)^l
f_1(\l)^{l-1}\;\;f_2(\l)^{r-l}\nonumber\\
&=O(\sqrt{n})\left(\frac{1}{2m - 2(l+r) }\right)^{l+r-1}
\left(\frac{er}{l}\right)^l \frac{(2c\lambda)^{r}}{\l f_2(\l)}
\left(\frac{2c\lambda f_1(\lambda)}{f_2(\l)^2}\right)^{l-1}\label{zx3}\\
&\text{using $\frac{\lambda f_2(\l)}{f_3(\l)} = 2c$}\nonumber\\
&\leq O(\sqrt{n})\left(\frac{1}{2m - 2(l+r) }\right)^{l+r-1}
\frac{2e(1-\b)c(2c\lambda)^{r}}{\l f_2(\l)}\left(\frac{4(1-\b)ec^2\lambda f_1(\l)}{f_2(\l)^2}\right)^{l-1}
\label{zx4}\text{\;\;\;\;using
$r/l \leq2(1-\b)c$}\\
&\leq O(\sqrt{n})\left(\frac{1}{2m - 2(l+r) }\right)^{l+r-1}
\frac{2e(1-\b)c(2c\lambda)^{r}}{\l f_2(\l)}\left(\frac{4(1-\b^2)ec^2\lambda}{f_2(\l)}\right)^{l-1}\text{\;\;\;\;using
$\frac{f_1(\l)}{f_2(\l)} < 1+\b$}\nonumber\\
&\leq O(\sqrt{n})\left(\frac{1}{2cn}\right)^{l+r-1}e^{3(l+r)^2/2cn}
(2c\lambda)^{r} \left(\frac{4ec^2\lambda}{f_2(\l)}\right)^{l-1}\text{\;\;\;\;using
$m = cn$}\nonumber\\
&=O(\sqrt{n})\left(\frac{1}{n}\right)^{l+r-1}e^{o(l)}
\;\lambda^{r} \left(\frac{2ec\lambda}{f_2(\lambda)}\right)^{l-1} \label{forgot}
\end{align}
since $r=O(l)$.

We now count the number of such configurations. We begin by choosing
\textit{EVEN(T)} and the root vertex of the tree in at most
$n\binom{n}{l-1}$ ways.  We make the following observation about
$T$. The contraction of
the lost edges of the tree yields a unique
tree on $l$ vertices. We note, by Cayley's formula, that the number of trees that
could be formed using $l$ vertices is $l^{l-2}$.  Reversing this contraction, we now
choose the sequence of $l$ vertices, \textit{Odd(T)}, that connect up
vertices in \textit{EVEN(T)} in $(n-l)(n-l-1)...(n-2l+1) =
(n-l)_l$ ways. We pick the remaining $r-l$ vertices from the
remaining $n-2l$ vertices in $\binom{n-2l}{r-l}$ ways. These
$r-l$ vertices can connect to any of \textit{EVEN(T)} in $l^{r-l}$
ways.  Hence, the total number of choices for $T$ is at most
\beq{zx2}
\binom{n}{l}l^{l-2} (n-l)_l \binom{n-2l}{r-l}l^{r-l}
\leq n^{r+l}e^{r} \left(\frac{l}{r-l}\right)^{r-l}.
\eeq

Combining the bounds for probability and choices of $T$, we get an upper bound of
\beq{zx5}
n^{r+l}e^r\left(\frac{l}{r-l}\right)^{r-l}
O(\sqrt{n})\left(\frac{1}{n}\right)^{l+r-1}e^{o(l)} \;\lambda^{r}
\left(\frac{2ec\lambda}{f_2(\l)}\right)^{l-1} \leq O(n^{3/2}) \cdot \left(\frac{ e\lambda l}{r-l}\right)^{r-l}
\left(\frac{4e^{2+o(1)}c\lambda^2 }{f_2(\l)}\right)^{l-1}
\eeq
The expression $\left(\frac{ e\lambda l}{x}\right)^{x}$ is maximized
at $x = \lambda l$. Our assumptions imply that $r\leq 2(1-\b)cl < \lambda l$.
Hence, we have the bound
\begin{align}
&O(n^{3/2}) \cdot \left(\frac{ e\lambda l}{2(1-\b)cl}\right)^{2(1-\b)cl}
\left(\frac{30c\lambda^2 }{f_2(\l)}\right)^{l}\nonumber\\
&\leq O(n^{3/2}) \cdot \left(\bfrac{ e }{1-\b}^{2(1-\b)c}\cdot
\frac{120c^3 }{f_2(\l)}\right)^{l}
\;\;\;\text{\;\;\; using
$\lambda < 2c $}\nonumber\\
&\leq O(n^{3/2}) \cdot e^{-\b cl}\label{zx9}
\end{align}
using 
$$f_2(\l) >\frac{ 120c^3e^{(2-\b)c}}{ (1-\b)^{2(1-\b)c}}$$
for $c$ sufficently large.

We sum $O(n^{3/2}) \cdot e^{-\b cl}$ over all $r$ and $l$ with
$L_0\leq l \leq n^{0.6}$ and $\; l \leq r
\leq (1-\b)cl$ and we get the probability to be at most

\beq{zx1}
O(n^{7/2}) e^{-\b cL_0} =o(1)
\eeq
for $c$ sufficiently large.

We now consider the probability of
the existence of a $T$ having $L_0 \leq l\leq n ^{0.6}$ and 
$r \geq 2(1+\b)cl$. Note that we can assume $r\leq l\D\leq l\log n$ here.

The bound \eqref{zx3} remains valid. Replacing $r$ by $r+1$ multiplies this by a factor $O(cn^{-1}e^{l/r})$ and so 
for this bound
we can just assume that $r=2(1+\b)cl$. This changes the $1-\b$ in \eqref{zx4} to $1+\b$ and we replace \eqref{forgot} 
by
$$O(\sqrt{n})\left(\frac{1}{n}\right)^{l+r-1}e^{o(l)}
\;\lambda^{r} \left(\frac{2ec(1+\b)^2\lambda}{f_2(\lambda)}\right)^{l-1}.$$
We re-use \eqref{zx2} and replace \eqref{zx5} by
\beq{today}
O(n^{3/2}) \cdot \left(\frac{ e\lambda l}{r-l}\right)^{r-l}
\left(\frac{4e^{2+o(1)}(1+\b)^2c\lambda^2 }{f_2(\l)}\right)^{l-1}.
\eeq
Our assumptions imply that $r=2(1+\b)cl > \lambda l$.
Hence, we have the bound
\begin{align*}
&O(n^{3/2}) \cdot \left(\frac{ e\lambda l}{2(1+\b)cl}\right)^{2(1+\b)cl}
\left(\frac{30(1+\b)^2c\lambda^2 }{f_2(\l)}\right)^{l}\\
&\leq O(n^{3/2}) \cdot \left(\bfrac{ e }{2(1+\b)}^{2(1+\b)c}\cdot
\frac{120(1+\b)^2c^3 }{f_2(\l)}\right)^{l}
\;\;\;\text{\;\;\; using
$\lambda < c $}\\
&\leq O(n^{3/2}) \cdot e^{-\b cl}
\end{align*}
using 
$$ f_2(\l) > \frac{120(1+\b)^2c^3e^{(2+4\b)c}}{(2(1+\b))^{2(1+\b)c}} $$
for $c$ sufficiently large.

We sum $O(n^{3/2}) \cdot e^{-\b cl}$ over all $r$ and $l$ with
$L_0\leq l \leq n^{0.6}$ and $r\geq 2(1+\b)cl$ and we get the probability to be at most

\beq{zx6}
O(n^{7/2}) e^{-\b cL_0} =o(1)
\eeq
for $c$ sufficiently large.

We next consider the case where $l\geq n^\e$ and we can only use $O(n^\e)$ late edges. 
We will use \eqref{zx2}, which is still an upper bound and only modify \eqref{expression}.
Let
$$b(d,d',1-\a)=\binom{d}{d'}\frac{((1-\a)m)_{d'}(\a m)_{d-d'}}{(m)_d}= \binom{d}{d'}(1-\a)^{d'}\a^{d-d'}
\brac{1+O\bfrac{\log^2n}{n}}$$
for $d\leq\D\leq \log n$.

We replace \eqref{expression} by
\begin{align}
&O(\sqrt{n})\left(\frac{1}{2m - 2(l+r) }\right)^{l+r-1}\times\nonumber\\
&\sum_{d_i\geq 3}\sum_{\substack {d_i' \leq d_i,\,i\in[l]\\
\sum_{i=1}^{l} d_i' = r+l-1}}
\left(\prod_{i=1}^{l} \frac{\lambda^{d_i}d_i'!}{d_i!f_3(\l)}b(d_i,d_i',1-\a)
 \prod_{i=l+1}^{2l-1}
\frac{\lambda^{d_i}d_i(d_i-1)}{d_i!f_3(\l)}\prod_{i=2l}^{r+l-1}
\frac{\lambda^{d_i}d_i}{d_i!f_3(\l)}(1-\a)\right) \label{zx10}\\
&= O(\sqrt{n})\left(\frac{1}{2m - 2(l+r) }\right)^{l+r-1}
\times\frac{\lambda^{2r+2l-3}(1-\a)^{2r-1}}{f_3(\lambda)^{r+l-1}}\times\nonumber\\
&\sum_{d_i\geq 3}\sum_{\substack{d_i'\leq d_i\\ \sum_{i=1}^{l} d_i' = r+l-1}}
\prod_{i=1}^l\frac{(\l\a)^{d_i-d_i'}}{(d_i-d_i')!}
\brac{\brac{\prod_{i=l+1}^{2l-1}\sum_{d_i\geq 3}\frac{\lambda^{d_i-2}}{(d_i-2)!}}
\brac{\prod_{i=2l}^{r+l-1}\sum_{d_i\geq 3}\frac{\lambda^{d_i-1}}{(d_i-1)!}}}\nonumber\\
&\leq O(\sqrt{n})\left(\frac{1}{2m - 2(l+r) }\right)^{l+r-1}
\times\frac{\lambda^{2r+2l-3}(1-\a)^{2r-1}}{f_3(\lambda)^{r+l-1}}\times\nonumber\\
&\sum_{\substack{d_i'\geq 1\\ \sum_{i=1}^{l} d_i' = r+l-1}}
\brac{\sum_{\k\geq 0}\frac{(\l\a)^{k}}{k!}}^l
\brac{\prod_{i=l}^{2l-1}\sum_{d_i\geq 3}\frac{\lambda^{d_i-2}}{(d_i-2)!}}
\brac{\prod_{i=2l}^{r+l-1}\sum_{d_i\geq 3}\frac{\lambda^{d_i-1}}{(d_i-1)!}}\nonumber\\
&\leq O(\sqrt{n})\left(\frac{1}{2m - 2(l+r) }\right)^{l+r-1}
\times\frac{\lambda^{2r+2l-3}(1-\a)^{2r}}{f_3(\lambda)^{r+l-1}}\times\binom{r}{l}e^{\l\a l}f_1(\l)^{l-1}f_2(\l)^{r-l}.
\label{zx8}
\end{align}
{\bf Explanation for \eqref{zx10}:} $d_i$ is the degree of vertex $i$ and for $i\in [l]$, $d_i'$ is the 
``early`` degree. The factor $b(d_i,d_i',1-\a)$ is the probability that $i$ has $d_i'$ neighbors.
 
Observe now that the expression in \eqref{zx8} is precisely 
$$e^{\l\a l}(1-\a)^{2r}\leq e^{\a(\l l-2r)}$$ 
times the expression in \eqref{zx7}. It follows that the probability bound \eqref{zx9} can be replaced by
$$O(n^{1/2}) \cdot e^{-\b cl}\cdot e^{\a(\l l-2r)}\leq O(n^{1/2}) \cdot e^{-\b cl/2}.$$ 
We sum this over $l,r$ to get the required conclusion.

The case $r\geq 2(1+\b)cl$ for $l\leq n^\e$, using only punctual edges follows a fortiori from the previous analysis.

We finally consider the case where $T$ is not a tree. 
When this happens, it will be because of a unique (Lemma \ref{lemgrow}) edge introduced in Case 1.
We can be handle this  by multiplying our final estimates 
by $O(L_0^2n^{-1}\log^2n)$. The factor $O(L_0^2)$ accounts for choosing a pair of vertices in $T$ in Case 1 
and $O(n^{-1}\log^2n)$ bounds the probability of the existence of this edge, given previous edges.

Part (c) follows from (b).

(d) If we consider the growth of the sub-tree emanating from $x$ then we can argue that it grows as fast as described in
(a) and (b). We just have to deal with the edges pointing into the part of $T$ that has already been constructed.
We can argue as for \eqref{zx10} with $\a=o(1)$, since the chances of choosing an endpoint in $T$ is $o(1)$ at each point.

If $x\notin W_1$ then the descendants $D_i$  of $x$ at levels $k_0+i$ grow at a rate of at least two 
(i.e. $|D_{i+1}|\geq 2|D_i|$) for $O(\log\log n)$ steps until $|D_i|\gg\log n$ and after this will grow
at a rate of al least $2c(1-\b)$. In which case the leaves of $T_x$, the sub-tree of $T$ rooted at $x$, will
constitute a fraction $1-O(1/c)$ of the vertices of $T_x$. The result now follows from Lemma \ref{lemgrow}(c). 

If $x\in W_1$ then $|D_i|$ grows at a rate of at most $2c(1+\b)$ once it has reached size $\log n$.
\proofend
\begin{remark}\label{rem2}
It follows from this lemma that only $O(n^{1/2+O(\e)})$ tardy edges are needed to build all
of the instances of $A_{k_0}$ needed by \HAM. If one looks at Section 4.3.1 of \cite{CFM} one sees, in conjunction with
equation (1) of that paper that the total running time of Step 3b of this paper is $O(n^{.995+o(1)})$ and so we can use
this as a bound on the number of punctual edges examined by Step 3b. We can drastically reduce this 
in the same way we did for building the trees in \HAM, but since we are only claiming our result for $c$ sufficiently 
large
and $\e\ll .005$, this is not necessary,
since there will \whp\ be $\Omega(n^{1-2\e})$ tardy $R_0:\wL$ edges, see Lemma \ref{tardy} below. 
In other words, almost all of the tardy $R_0:\wL$ edges are {\bf not} using for tree 
building. 
\end{remark}

The above lemma shows that $A_k$ can be relied on to get large. Unfortunately, we need to do some more analysis because 
we do not have full independence, having run \2G. Normally, one would only have to show that $END(a)$ is large for
all relevant vertices $a$ and this would be enough to show the existence \whp\ of an edge joining $a$ to $b\in END(a)$ 
for some $a,b$. We will have to restrict our attention to the case where $a\in R_0$ and $b\in \wL$, see Remark
\ref{rem1}. So first of all we will show that \whp\ there are many $a\in R_0$, see Lemma \ref{lem2}. 
For this we need to show that every path we come across contains many consecutive triples $u,v,w\in R_0$. In which
case, an inserted edge $(x,v)$ produces a path with an endpoint in $R_0$. 
We also need to show 
that \whp\ there are many $b\in \wL$, see Lemma \ref{lem:tree-not-all-early}. We will also need to show that there are
many edges that can be $(a,b)$, see Lemma \ref{tardy}.

For the Lemma \ref{conn1} below we need some results from \cite{F2}.
Let $\bu=\bu(t)$ denote $(y(t),z(t),\m(t))$ and let $\hu=\hu(t)$ denote $(\hy(t),\hz(t),\hm(t))$ where $y(t)$ etc. denotes
the value of $y=|Y|,z=|Z|,\m=|E(\G(t)|$ at time $t$ and 
$\hy(t)$ etc. denotes the deterministic value for the solution to the associated
set of differential equations, summarised in equation (152) of that paper:
\beq{slide}
\frac{d\hy}{dt}=\hA+\hB-\hC-1;\quad\frac{d\hz}{dt}=
2\hC-2\hA-2\hB;\quad\frac{d\hm}{dt}=-1-\hD.
\eeq
where
\beq{ABCD}
\hA=\frac{\hy\hz\hl^5f_0(\hl)}{8\hm^2f_2(\hl)f_3(\hl)},\quad \hB=\frac{\hz^2\hl^4f_0(\hl)}
{4\hm^2f_2(\hl)^2},
\quad \hC=\frac{\hy\hl f_2(\hl)}{2\hm f_3(\hl)},\quad \hD=\frac{\hz\hl^2f_0(\hl)}{2\hm f_2(\hl)}.
\eeq
and $f_j(x)=e^x-\sum_{i=0}^{j-1}\frac{x^i}{i!}$.

Lemma 7.1 of \cite{F2} proves that $\bu(t)$ and $\hu(t)$ are close \whp:
\begin{lemma}\label{close}
$$||\bu(t)-\hu(t)||_1\leq  n^{8/9}, \qquad\text{ for }1\leq t\leq \min\set{T_0,\hT_0}\ \whp.$$
\end{lemma}
Here $T_0$ is a stopping time and $\hT_0$ is a deterministic time such that \whp\ Step 3 begins before 
$\min\set{T_0,\hT_0}$.

Note that $\e\ll1/9$. 
Let 
$$i_o=n^{3/4-\e}\text{ and }\r=n^{1/4}.$$
Equation (163) of \cite{F2} states that \whp
\beq{xi}
|\th_\xi(\hu(t))-\D_\xi|=O\brac{\r^{-1}\log^2n+\frac{||\bu(t)-\hu(t)||_1}{n}}
\text{ for }\xi=a,b,c,2.
\eeq
Here
$$\th_a=0,\;\th_b=\hA\;\th_c=\hA+\hB\text{ and }\th_2=1-\th_a-\th_b-\th_c$$
and $\D_\xi$ is the proportion of steps in $[t,t+\r]$ that are Step 1$\xi$ or Step 2, if $\xi=2$.

Now if $\hz=o(n)$ $\hm=\Omega(n)$and $\hl=\Omega(1)$ then we have from \eqref{ABCD} that $\hA,\hB=O(\hz/n)$ and that 
$\th_b=O(\hz/n),\th_c=O(\hz/n),\th_2=2-o(1)$.
Then from \eqref{slide} we see that $\hz$ grows at the rate $2-o(1)$ per time step, so long as $t=o(n)$ and
hence $\hz=o(n)$. 

It is shown in \cite{F2} that if $c\geq 10$
then \whp\  $\hl=\Omega(1)$ up until the (random) time when Step 3 begins. 
See equation (190) of that paper. Furthermore, it follows from Lemma
\ref{close} that \whp
\begin{enumerate}[{\bf X1}]
\item If $t=\g n^{1-\e}$ for some constant $\g$ then \whp\ $z(t)\sim 2t$.
\item If $t=\g n^{1-\e}$ for some constant $\g$ then \whp\ there will be $O(n^{1-2\e})$ instances of Step 1 in 
$[0,t]$.
\item $\l=\Omega(1)$ up until the start of Step 3.
\end{enumerate}

\begin{lemma}\label{conn1}
W.h.p., all the paths in Steps 1 and 2 of \HAM\ contain at least $n_0=\Omega(n^{1-4\e}/\log n)$ pairs of consecutive 
edges $(u,v),(v,w)$ such that $u,v,w\in R_0$.
\end{lemma}
\proofstart
First consider the steps in the range $[0,i_0\r/4]$. It follows from 
{\bf X1} that at the end of this period, there will \whp\ be at least $i_0\r/3$
vertices in $Z$. 
Consider the edge $(v,w)$ of Step 2 at some time in $[0,i_0\r/4]$. The probability that $w\in Z$ is certainly 
$\Omega(n^{-\e})$ and the probability it has a punctual $Z$-witness is $1-\a-o(1)$. This holds regardless of the 
previous history,
once we condition on an event that happens \whp

The probability that $w\in Z$ is certainly 
$O(n^{-\e}\D)=O(n^{-\e}\log n)$.
This implies that the number of times we create a component of $M$ containing more than two vertices is
$O(n^{1-2\e}\log n)$.
Thus almost all components of $M$ at the  
end of the period $[0.i_0\r/4]$ consist of isolated edges. Let us assume then that there are at least $A_1n^{1-\e}$ such 
edges
where in the following $A_1,A_2,\ldots,$ are positve constants. Let $S_1$ denote this set of components.

Now consider the steps in the range $[i_0\r/4,i_0\r/2]$ and consider the edge $(v,w)$ of Step 2. 
We have $w\in V(S_1)$ with probability at least $A_2 n^{-\e}$. This is because \whp\ the 
total degree of $V(S_1)$ will be
$\Omega(n^{1-\e})$ and the total degree of $G$ is at most $2cn$. The vertex $w$ is early by construction. Also the 
$Z$-witness of $w$ will be punctual
with probability at least $1-\a-O(n^{-\e})$. We next observe that with probability at least 
$\brac{1-\Omega\bfrac{\D}{n}}^{n^{1-\e}/4}=1-o(1)$, this component will not be absorbed into a 
larger component in $[i_0\r/4,i_0\r/2]$. 
Thus, in expectation, at time $i_0\r/2$ there is a set of $A_3n^{1-2\e}$ components of $M$ consisting of a path
of length two with its middle vertex in $R_0$. A simple second moment calculation will show concentration around the
mean, for $|S_2|$. 

We can repeat this argument for the periods $[i_0\r/2,3i_0\r/4]$,$[3i_0\r/4,i_0\r]$ to argue that by time 
$n^{1-\e}$, $M$ will contain a set $S_3$ of at least $A_4n^{1-4\e}$ components consisting of paths of length 
four in which
the internal vertices are all in $R_0$.

We can argue that \whp\ at least half of the components in $S_3$ will have both end vertices of degree at most $3c$.
Denote these by $S_4$.
Indeed the number of edges incident with vertices of degree more than $3c$ is relatively small.
Indeed, the expected number of such edges is asymptotically equal to $\sum_{k\geq 3c}\frac{k\l^k}{k!f_3(\l)}\leq
\e_c=(e/3)^{3c}$. The number of such edges is concentrated around its mean. If we assume degrees are independent and 
less than $\log n$ then
we can use Hoeffding's Theorem and then correct by a factor $O(n^{1/2})$ to condition on the total degree. Given this,
we see that \whp\ at least $2(1-\e_c)^2/3$ of the components of $S_3$ will be created in two executions of Step 2 with
the degree $v$ less than $3c$. 

Observe now that with probability at least $\brac{1-\frac{6c}{\Omega(n)}}^{2cn}=\Omega(1)$ a component $C\in S_3$
will survive as a component of $M$ until the end of Step 2. 
Because $|S_3|$ is $O(n^{1-\e})$, this is true regardless of which other compponents in $S_3$ survive.
The $\Omega(n)$ in the denominator comes from the fact that
\whp\ Step 2 ends with $|Z|=\Omega(n)$. Let $S_4$ denote this set of components and note that \whp\ there will be
at least $A_5n^{1-4\e}$ components in $S_4$. 

Step 3 of \2G\ adds a matching $M^{**}$ that is disjoint from the edges in the contraction of $S_3$ to a matching.
This matching is independent of $S_3$. This implies that \whp\ any cycle (or possibly path) of the union of 
$M^*$ and $M^{**}$ of length $\ell\geq n^{8\e}$, contains at least $A_5\ell n^{-4\e}$ members of $S_4$. Here we are
using concentration of the hypergeometric distribution i.e. sampling without replacement.

In \HAM\ we start with a path of length $\ell=\Omega(n/\log n)$ and \whp\ every path is generated by deleting at
most $O(\log^2n)$ edges. This completes the proof of the lemma.
\proofend

\subsection{Batches}
Let $\G(t)$ denote the graph $\G$ after $t$ steps of \2G.
Suppose that $t_1<t_2\leq n^{1-\e}$ and that \2G\ applies Step 2 at times $t_1,t_2$ and Step 1 at times $t_1<t<t_2$.
We consider the
set of edges and vertices removed from time $t_1$ to time $t_2$, i.e.  the graph
$\G(t_1)\setminus \G(t_2)$ and call it a {\em batch}. Note that batches are connected subgraphs
since each edge/vertex removed is incident to some edge that is also removed. 

We also claim that each batch
\whp\ is constructed within $O(\log^2n)$ steps and 
contains $O(\log^3n)$ vertices. This follows from \cite{F2} as we now explain. 
Let $\z=y_1+2y_2+z_1$. Equations (67), (68), (69) of \cite{F2} show that 
$$\ex[\z'-\z\mid|\bv|]=-(1-Q)-o(1)$$
where 
$$Q=Q(\bv)=\frac{yz}{4\mu^2}\,\frac{\la^3}{f_3(\la)}\,\frac{\la^2f_0(\la)}{f_2(\la)}+
\frac{z^2}{4\mu^2}\,\frac{\la^4f_0(\la)}{f_2(\la)^2}.$$
Lemma 6.2 of \cite{F2} shows that $1-Q=-\Omega(1)$ if $\l=\Omega(1)$, and {\bf X3} is our justification for assuming this. 

Thus the expected change in $\z$ is $-\Omega(1)$ when $\z>0$. 
We carry out Step 2 iff $\z=0$. Now $\z$ can change by at most $O(\D)=O(\log n)$ and has a negative drift whenever
it is positive. This implies that it
must return to zero within $O(\log^2n)$ steps. 
Another $\D\leq \log n$ factor will allow at most $\log n$ edges to be removed in one step.
By making the hidden constant sufficently large, we can replace \whp\ by with probability $1-O(n^{-10})$.

\begin{lemma}\label{lem:batch-graph-edges}\ 
\begin{enumerate}[(a)]
\item W.h.p. there are at most $n^{1-4\e}$ vertices $v\in G$ that are within
distance $\ell_0=2\log\log n$ of
6 distinct batches.
\item W.h.p. no vertex has degree more than 4 in a single batch.
\end{enumerate}
\end{lemma}
\proofstart\ 

(a) We bound the probability of being within distance $\ell_0$ of $s$ batches by
$$\r_s=\sum_{v=1}^n\binom{n^{1-\e}}{s}\prod_{i=1}^{s}\Pr(dist(v,B_i)\leq \ell_0
\mid dist(v,B_j)\leq \ell_0,\,1\leq j<i).$$
{\bf Explanation:} Here $\binom{n^{1-\e}}{s}$ is the number of choices
for the start times of the batches\\ $B_1,B_2,\ldots,B_s$.

We claim that for each $i,v$,
\beq{swim}
\Pr(dist(v,B_i)\leq \ell_0\mid dist(v,B_j)\leq \ell_0,\,1\leq j<i)=
O\bfrac{\log^{2+\ell_0}n}{n}.
\eeq
This gives
$$\r_s\leq \exp\set{-(K-2-o(1))(\log\log n)^2s}.$$
and assuming $K\geq 13$,
this implies that the expected number of vertices within distance $\ell_0$ of 6 batches is less than $n^{1-5\e}$.
The result now follows from the Markov inequality.

{\bf Proof of \eqref{swim}:} Suppose that $B_i$ is constructed at time $t_i$. It is a subgraph of $\G(t_i)$
and depends only on this graph. We argue that
\beq{atlast}
\Pr(\exists w\in B_i:\,dist(v,w)\leq \ell_0\mid dist(v,B_j)\leq \ell_0,\,1\leq j<i)
\leq O(n^{-10})+O\bfrac{i\log^{2+\ell_0}n}{n}.
\eeq
{\bf Explanation:} The $O(n^{-10})$ term is the probability the batch $B_i$ is large.
The term $O\bfrac{i\log^{2+\ell_0}n}{n}$ in \eqref{atlast} arises as follows. We can assume that
$|N_{\ell}(v)|\leq \D^{\ell_0}\leq\log^{\ell_0} n$, where $N_\ell(v)$ is the set of vertices within distance
$\ell$ of $v$. Suppose as in \cite{AFP} we expose the
graph $\G$ at the same time that we run \2G. For us it is convenient to work within the configuration model of
Bollob\'as \cite{B2}. 
Assume that we have exposed $N_\ell(v)$. At the start of the construction of a batch we choose a
random edge of the current graph. The probability this edge lies in $N_{\ell_0}(w)$ is
$O(\log^{\ell_0} n/n)$. In the middle of the construction of a batch,
one endpoint of an edge is known and the the other endpoint is chosen randomly from the set of configuration
points associated with $\G(t)$. The probability this new endpoint lies in $N_{\ell_0}(v)$ is
also $O(\log^{\ell_0} n/n)$
and there are only $O(\log^2n)$ steps in the creation of a batch. 


(b) The probability that vertex $v$ appears $k+3$ times in a fixed batch can be bounded above by\\
$\binom{O(\log^2n)}{k}\brac{O\bfrac{\D}{n}}^{k}=O\bfrac{\log^{k+3}n}{n^k}$. 
Indeed, if $v$ has degree at least 3 at any time, then the probability
its degree in the current batch increases in any step is $O\bfrac{\D}{n}$. 
\proofend

We now argue that there will be a sufficient number of tardy $R_0:\wL$ edges.
\begin{lemma}\label{tardy}
W.h.p. there will be $\Omega(n^{1-2\e})$ tardy $R_0:\wL$ edges.
\end{lemma}
\proofstart
We first consider the set $F_1$ of tardy edges $e=(u,v)$ such that 
(i) $u$ appears at least twice in the first $n^{1-\e}/10$ edges and in at least 30 other punctual edges and
(ii) vertex $v$ has degree at least 30 and does not appear in the first $n^{1-\e/2}$ edges in $\s$. 
It is straightforward to show 
that \qs\
we have $|F_1|=\Theta(mn^{-2\e})$. 

Suppose that $u$ satisfies (i). It loses at most 24 edges (Lemma \ref{lem:batch-graph-edges}(a),(b))
before the second edge incident with $u$ is chosen  and
then $u$ will be in $R_0$. This is because $u$ will be in $Z$ just before this point and will then be placed in 
$R$. and it will have at least six choices for a punctual $Z$-witness.
We use the fact that almost all of the first $n^{1-\e}$ steps are Step 2 to see that the edges incident with $u$
occuring in the first $n^{1-\e}/10$ steps will indeed be selected before time $n^{1-\e}$.
\ignore{
On the other
hand, there will \whp\ be $O(n^{1-4\e}\log^2n)$ vertices that lose two edges this way. Indeed, the probability of losing
$t$ edges this way is $O\brac{\binom{n^{1-2\e}\D}{t}n^{-t}}$. Thus \whp\ at most $O(n^{1-4\e}\log^2n)$ edges
of $F_1$ are not in $R_0:\wL$ because of this. Recall that \qs\ there are only $O(n^{1-2\e})$ executions of Step 1
in the first $n^{1-\e}$ steps of \2G, see {\bf X2}. Note that the vertices that satisfy (i) and lose more than 24 edges
before the second edge incident with $u$ is chosen account for only $O(n^{1-4\e}\D)=o(|F_1|)$ pairs in $F_1$.
This follows from Lemma \ref{lem:batch-graph-edges}.}

If $v$ satisfies (ii) and loses at most 24 edges because of Step 1 in the first $n^{1-\e}$ steps
then $v$ will be in $\wL$. This is because it will have degree at least six in $\G_{n^{1-\e}}$.
\ignore{
The number of vertices satisfying (ii) and losing 24 edges because of Step 1 
accounts for $O(n^{1-4\e}\D)=o(|F_1|)$ pairs in $F_1$.}
\proofend

We now consider the probability that $A_{k_0}$ contains many vertices that lie in $\Lambda_1=\Lambda_2\cup
\Lambda_3$ where 
\begin{align*}
&\Lambda_2=\set{v: \text{$v$ appears in the first $n^{1-\e/2}$ edges in $\s$}}\\
&\Lambda_3=\set{v\notin \Lambda_2: \text{$v$ loses 24 edges because of Step 1 in the first $n^{1-\e}$ steps}}.
\end{align*}
In the proof of Lemma \ref{tardy} we used the fact that 
\beq{usef}
\text{if the degree of $v$ is at least 30 and }v\notin \Lambda_1\text{ then }v\in\wL.
\eeq
\begin{lemma}\label{f**k}
W.h.p., every extension-rotation tree $T$ has $|A_{k_0}\cap \Lambda_1|\leq |A_{k_0}|/30$.
\end{lemma}
\proofstart
We first estimate $|A_{k_0}\cap \Lambda_2|$. We go back to \eqref{today} and estimate, for fixed $r,l$, 
\beq{cold}
\Pr(\exists T: |A_{k_0}\cap \Lambda_2|\geq l/100)\leq O(n^{3/2}) \cdot \left(\frac{ e\lambda l}{r-l}\right)^{r-l}
\left(\frac{4e^{2+o(1)}(1+\b)^2c\lambda^2 }{f_2(\l)}\right)^{l-1}\binom{l}{l/200}\bfrac{n^{1-\e/2}\D}{m}^{l/200}.
\eeq
{\bf Explanation:} 
\ignore{We have ignored a term $o(n^{-10})$ say, that accounts for maximum degree $\leq\log n$. Given this,}
We have taken the RHS of \eqref{today} and multiplied by a bound on the 
probability that there are at least $l/200$ members of
$EVEN(T)$ appearing in the first $n^{1-\e/2}$ edges of $\s$. Note the the permutation $\s$ is independent of $T$ and
that \whp\ we will have $|A_{k_0}|\geq l/2$. This is because for most of the time, the tree grow at a rate at
least $2c(1-\b)$. We use \eqref{today} and not \eqref{zx5} because we can only assume that $r\leq 2(1+\b)cl$.

Thus,
\begin{multline}\label{ffdd}
\Pr(\exists T: |A_{k_0}\cap \Lambda_2|\geq l/100)\leq \\O(n^{3/2})\cdot\brac{e^{\l l/(r-l)}\cdot 
\frac{4e^{2+o(1)}(1+\b)^2c\lambda^2 }{f_2(\l)}\cdot (200e)^{1/200}\cdot\bfrac{n^{1-\e/2}\log n}{cn}^{1/200}}^l
\leq n^{-\e l/300}.
\end{multline}
Now we are interested here in the case where $l=n^{1/2+o(1)}$ and so this is easily strong enough so that we can apply
the union bound over $r,l$.

We next estimate $|A_{k_0}\cap \Lambda_3|$. We replace \eqref{cold} by
\begin{multline}\label{colder}
\Pr(\exists T: |A_{k_0}\cap \Lambda_3|\geq l/100)\leq \\
O(n^{3/2}) \cdot \left(\frac{ e\lambda l}{r-l}\right)^{r-l}
\left(\frac{4e^{2+o(1)}(1+\b)^2c\lambda^2 }{f_2(\l)}\right)^{l-1}\binom{l}{l/200}\binom{n^{1-\e}}{23l/200}
\bfrac{l\log n}{200m}^{23l/200}.
\end{multline}
{\bf Explanation;} We have taken the RHS of \eqref{today} and multiplied by a bound on the 
probability that there is a set of leaves $S$ of size $l/200$ such that at least $24l/200$ times during the 
first $n^{1-\e}$ steps the vertex $w$, the neighbor of the selected $v$, is in $S$. 
This is computed assuming that we have exposed the edges of $T$. Note that we have that in at most $l/200$ times
do we lose the edge of $T$ incident with $w\in S$, explaining the factor $\binom{n^{1-\e}}{23l/200}$.

Equation \eqref{ffdd} is replaced by
\begin{multline*}
\Pr(\exists T: |A_{k_0}\cap \Lambda_3|\geq l/100)\leq\\ 
O(n^{3/2}\cdot\brac{e^{\l l/(r-l)}\cdot 
\frac{4e^{2+o(1)}(1+\b)^2c\lambda^2 }{f_2(\l)}\cdot (200e)^{1/200}\cdot\bfrac{200en^{1-\e}}{23l}^{23/200}\cdot
\bfrac{l\log n}{200cn}^{23/200}}^l\leq  n^{-\e l/20}.
\end{multline*}
\proofend

The next lemma puts a lower bound on $|\wL\cap A_{k_0}|$ (see Lemma \ref{lem1}).

\begin{lemma}\label{lem:tree-not-all-early}
W.h.p. $|\wL\cap A_{k_0}|\geq |A_{k_0}|/2$.
\end{lemma}
\proofstart
Let $k_1=k_0-2\ell_0$ where $\ell_0=2\log\log n$. 
and consider $A_{k_1}=\set{a_1,a_2,\ldots,a_r}$. Note that 
$$r\geq \frac{n^{1/2+5\e}}{2c(1+\b)\log^{\ell_0}n}.$$
Let $s_i$ be the number
of descendents of $a_i$ in $A_{k_0}$ and let $s_i'$ be the number of
early descendents of $a_i$ in $A_{k_0}\cap \Lambda_1$.

Let $s_i''$ be the number of descendents of $a_i$ in $A_{k_0}$ that
have degree at most $30$. 
We observe from Lemma \ref{lem1}(b) that
\beq{ll1}
|A_{k_0}|=\sum_{i=1}^rs_i\geq r(2c(1-\b))^{k_0-k_1}\geq r\log^{\log c}n.
\eeq
Next let $I=\set{i\in [r]: a_i\notin W_1\text{ and }s_i\geq (2c(1-\b))^{k_0-k_1}/4}$ (where $W_1$ is from
Lemmas \ref{lemgrow}, \ref{lem1}) and observe that
\beq{qazx}
\sum_{i\notin I}s_i\leq r(2c(1-\b))^{k_0-k_1}/4+n^{1/2}\log^{4\ell_0+1}n(2c(1+\b))^{k_0-k_1}\leq |A_{k_0}|/3.
\eeq
It follows from Lemma \ref{lem1}(d) that 
$$s_i''\leq s_i/5\qquad\text{ for }i\in I.$$

It follows from Lemma \ref{f**k} that \whp
$$\sum_{i=1}^rs_i'\leq|A_{k_0}|/30.$$
Now, after using \eqref{usef}, we see that
$$|\wL\cap A_{k_0}|\geq \sum_{i\in I}(s_i-s_i'')-\sum_{i=1}^rs_i'\geq \brac{\frac45\cdot\frac23-\frac{1}{30}}|A_{k_0}|.$$
\proofend

We now consider going one iteration further and building $A_{k_0+1}$.
\begin{lemma}\label{lem2}
W.h.p. $A_{k_0+1}$ contains at least $\Omega(n^{1/2+2\e})$ vertices of $R_0$.
Furthermore, we can find these $R_0$ vertices by examining $n^{1-3\e}\log n$ tardy $R_0:\wL$ edges.
\end{lemma}
\proofstart
Assume from Lemmas \ref{lem1} and \ref{lem:tree-not-all-early} that $A_{k_0}$ contains at least
$n_1=\frac{n^{1/2+5\e}}{4c^2(1+\b)}$ vertices in $\wL$. Assume also from Lemma \ref{conn1} that all of the 
paths corresponding to $A_{k_0}$ have $n_0=\Omega(n^{1-4\e}/\log n)$ 
consecutive triples $u,v,w\in R_0$. If the middle vertex
$v$ is the neighbour of an endpoint, then it yields a new endpoint of $A_{k_0+1}$ in $R_0$.
Then the expected number of rotations
leading to an endpoint in $R_0$ is at least
$$C_1\times n^{1-3\e}\log n\times \frac{n^{1/2+5\e}}{3c^2(1+\b)}\times \frac{n^{1-4\e}}{\log n}
\times \frac{1}{n^{1-2\e}\times n}=\Omega(n^{1/2+2\e})$$
for some constant $C_1>0$.

We can claim a \qs\ lower bound because almost all of the tardy $R_0:\wL$ edges are unconditioned, 
see remark \ref{rem2}.
\proofend
\section{Finishing the proof}\label{ftp}
We have argued that we only need to do 
$\ell_1=O(\log^2n)$ extensions w.h.p. The tardy $R_0:\wL$ edges are our scarce resource of residual
randomness. Remark \ref{rem2} explains that we only need to use an $o(1)$ proportion in building 
trees up to the $k_0$th level. We will only use the result of Lemma \ref{lem2} for growing the first
extension-rotation tree of each of the $O(\log^2n)$ path extensions. Lemma \ref{lem2} tells us that
we only need to use an $o(1)$ fraction of the available $R_0:\wL$ edges for producing many paths that have
an $R_0$ endpoint. 

Consider a round of \HAM\ where we are trying to extend path $P$. We start with a path and then 
we construct a BFS ``tree''. 
After the first tree construction of each round, we construct $A_{k_0}$ and 
create one more level $A_{k_0+1}$. From Lemma \ref{lem2}, we should obtain 
$\Omega(n^{1/2+2\e})$ paths with early
endpoints. Now we grow trees from each of these paths and try to close them using the 
set $E_L=\set{f_1,f_2,\ldots,f_M}$  of unused tardy $R_0:\wL$ edges. We can examine these edges in $\s$ order. 
The probability that
the next edge $f_i$ fails to close a path to a cycle is $p=\Omega(n^{1/2+2\e}\times n^{1/2+5\e}\times n^{-2})$.
So the probability we fail is at most $\Pr(Bin(M,p)<\ell_1)$. Now $Mp=\Omega(n^{5\e})\gg\ell_1$ and so the Chernoff
bounds imply that we succeed \whp

As final thought, although we have proved that we can find a Hamilton cycle quickly, 
being very selective in our choice of edges for certain purposes, the breadth first
nature of our searches imply that we can proceed in a more natural manner and use all 
edges available to us. In the worst-case we would have to use the designated ones.
\section{Why not $\e=O\bfrac{\log\log n}{\log n}$?}\label{diff}
In the proof of Lemma \ref{lem1} we need to choose $\ell_0=2\log\log n$ so that $2^{\ell_0}\gg L_0$ of that lemma.
But then in \eqref{swim} we want $n^{\e}\gg \log^{\ell_0}n$. With some work we could replace the bound $\log^{\ell_0}n$
by $O(c)^{\ell_0}$ which would allow us to take $\e=\frac{K\log\log n}{\log n}$. The catch here is that in this 
case we would need $K$ to grow with $c$. This is not satisfactory and so we content ourselves for now
with \eqref{defeps}.
\section{Final Remarks}\label{FM}
We have shown that a Hamilton cycle can \whp\ be found in $O(n^{1+o(1)}))$ time. 
It should be possible to replace $n^{o(1)}$ by $\log^{O(1)}n$ and we have explained the technical difiiculty
in Section \ref{diff}.
We think that $O(n\log^2n)$ should be possible. It should also be possible to apply the ideas here to speed up the known
algorithms for random regular graphs, or graphs with a fixed degree sequence.

\end{document}